\documentstyle[12pt]{article}
\begin{document}
\newpage
\pagestyle{empty}
\newfont{\twelvemsb}{msbm10 scaled\magstep1}
\newfont{\eightmsb}{msbm8}
\newfont{\sixmsb}{msbm6}
\newfam\msbfam
\textfont\msbfam=\twelvemsb
\scriptfont\msbfam=\eightmsb
\scriptscriptfont\msbfam=\sixmsb
\catcode`\@=11
\def\Bbb{\ifmmode\let\next\Bbb@\else
  \def\next{\errmessage{Use \string\Bbb\space only in math mode}}\fi\next}
\def\Bbb@#1{{\Bbb@@{#1}}}
\def\Bbb@@#1{\fam\msbfam#1}
\newfont{\twelvegoth}{eufm10 scaled\magstep1}
\newfont{\tengoth}{eufm10}
\newfont{\eightgoth}{eufm8}
\newfont{\sixgoth}{eufm6}
\newfam\gothfam
\textfont\gothfam=\twelvegoth
\scriptfont\gothfam=\eightgoth
\scriptscriptfont\gothfam=\sixgoth
\def\frak{\frak@}
\def\frak@#1{{\fam\gothfam{{#1}}}}
\def\frak@@#1{\fam\gothfam#1}
\catcode`@=12
\def\C{{\Bbb C}}
\def\N{{\Bbb N}}
\def\Q{{\Bbb Q}}
\def\R{{\Bbb R}}
\def\Z{{\Bbb Z}}
\def\id{\mbox{id}}
\newtheorem{theorem}{Theorem}[section]
\newtheorem{claim}[theorem]{Claim}
\newtheorem{proposition}[theorem]{Proposition}
\newtheorem{definition}[theorem]{Definition}
\newtheorem{lemma}[theorem]{Lemma}
\newtheorem{question}[theorem]{Question}
\newtheorem{corollary}[theorem]{Corollary}
\newtheorem{remark}[theorem]{Remark}
\newtheorem{example}[theorem]{Example}
\newtheorem{examples}[theorem]{Examples}

$$
\;
$$

\vskip 2cm

\begin{center}

  {\LARGE {\bf {\sf Abelian topological groups without
irreducible Banach representations}}} \\[1cm]

\smallskip

{{\large Vladimir Pestov} (Wellington, New Zealand)}
\vskip .6cm

To Professor Adalberto Orsatti's 60-th Birthday
\end{center}
\vskip .3cm

\begin{abstract}
We exhibit abelian topological 
groups admitting no nontrivial strongly continuous 
irreducible representations in Banach spaces. 
Among them are some abelian Banach--Lie groups 
and some monothetic subgroups of
the unitary group of a separable Hilbert space.
\bigskip

\bigskip

\end{abstract}

\pagestyle{plain}

\begin{section}{Introduction}
\mbox{}

The classical 1943 Gelfand--Ra\u\i kov theorem
\cite{GR} states that every
locally compact group has a complete system of irreducible strongly
continuous representations in Hilbert spaces. Does this result admit a
sensible generalization to wider classes of topological groups?

The most immediate problem is that
even the `nicest' non locally compact topological groups known to date
not necessarily possess nontrivial representations 
in Hilbert spaces, irreducible or not
 --- for example, there are abelian Banach--Lie groups without
weakly continuous Hilbert space representations \cite{Ba}. 
If one allows for
representations in more general Banach spaces, the situation becomes
more favourable. 
In 1957 Teleman \cite{T} proved that every Hausdorff topological
group admits a faithful strongly continuous representation 
in a suitable Banach space by isometries.
He further asked whether every topological
group admits a complete system of strongly continuous irreducible
representations in Banach spaces by isometries.
The answer is `No:' as noticed recently by
the present author \cite{Un}, no minimally almost periodic 
monothetic topological group
admits nontrivial irreducible representations in Banach spaces by
isometries. Can Teleman's conjecture survive if one drops the
restrictive
requirement `by isometries' and allows any strongly continuous Banach
representations?
The aim of this note is to show that the answer is still in the
negative. 

We exhibit a vast class of abelian topological 
groups admitting no infinite-dimensional continuous 
irreducible representations in Banach spaces: an abelian
topological group is such whenever the torsion subgroup is 
everywhere dense in it.
Clearly, every minimally almost periodic
abelian topological group (that is, one without nontrivial continuous 
characters)
contained in this class admits no
nontrivial irreducible Banach representations. We observe
that at least three previously known classes of examples of
minimally almost periodic groups fall into this category. 
Among them are all minimally almost periodic abelian Banach--Lie 
groups,
and even some monothetic subgroups of
the unitary group of a separable Hilbert space with the strong
operator topology, which fact shows that the
unitary representations of some of the 
`nicest' non locally compact groups known 
cannot in general be decomposed into irreducibles in any reasonable sense.

To substantiate the subject of this note, let us remark that
infinite-dimensional
irreducible Banach representation of the infinite cyclic group
do exist indeed: every bounded linear
operator on a Banach space having no invariant subspaces
\cite{E, Re, B} leads to such a representation. 
It remains yet to be seen if there
exists a minimally almost periodic
topological group admitting an irreducible continuous Banach
representation. The importance of this open problem (and especially
its version for Hilbert space representations) stems from
its obvious relevance to the Invariant Subspace Problem.
Our present results suggest, however, that learning to produce representations
of this kind might take more than a straightforward extension of the 
Gelfand--Ra\u\i kov theorem.
\end{section}

\begin{section}{Main result}
\mbox{}

A representation $\rho$ of a topological group $G$ in a normed space $E$
is {\it strongly continuous,} or simply
{\it continuous,} if it is continuous as a mapping 
$\rho\colon G\times E\to E$. Equivalently, $\rho$ is strongly
continuous if it is continuous as a homomorphism
$\rho\colon G\to\mbox{GL}(E)$, where the general linear group
$\mbox{GL}(E)$ of $E$ is equipped with the strong operator topology.
It is worth remembering
 that the strong operator topology is never a group topology
on $\mbox{GL}(E)$, but it is such on the subgroup formed by all
isometries, in particular, on the unitary group $U(\cal H)$ of
a Hilbert space $\cal H$.

Recall (see e.g. \cite{Ka})
that the collection of all elements of an abelian group
$G$ having finite order forms a subgroup of $G$ called the
{\it torsion subgroup,} which we denote by $T(G)$.
If every element of $G$ is of finite order,
then $G$ is called a {\it torsion group.}

\begin{theorem}
Let $G$ be a topological abelian group which is algebraically a
torsion group. Then every continuous irreducible
representation of $G$ in a Banach space is one-dimensional and indeed a
continuous character.
\label{t1}
\end{theorem}

\noindent{\bf Proof.} Let $\rho$ be a 
continuous representation of $G$ in a Banach space $E$.
Assume that for some $g\in G$, the 
operator $\rho_g$ is not a scalar multiple of the identity.
Now a standard argument from operator theory 
(see e.g. \cite{B}, p. 277) shows that $\rho$ admits a
proper subrepresentation.
Indeed, since $(\rho_g)^n=\mbox{I}$ for some $n$, the bounded linear
operator $\rho_g$ has an eigenvalue, $\lambda$, which is an $n$-th
root of unity. (See e.g. \cite{Ka}, Sect. 12, or \cite{B}, Exercise
2, p. 35.)
The space $F=\mbox{Ker}(\lambda\mbox{I}-\rho_g)$
is closed, proper (since $\rho_g\neq\lambda\mbox{I}$), and clearly
invariant under every operator from ${\cal L}(E)$ commuting with
$\rho_g$. (As operator theorists say, $F$ is
{\it hyperinvariant}.)
Since $G$ is abelian, the space $F$ is invariant under every
operator $\rho_{h}$, $h\in G$, and $\rho\vert_{F}$ forms a proper
Banach subprepresentation of $\rho$. 

We conclude that if $\rho$ is irreducible, then each $\rho_g$,
$g\in G$ is a scalar multiple of the identity and 
$\dim E=1$. Since $\rho(G)$ is a torsion multiplicative
subgroup of $\C^{\times}$,
it is contained in $\Bbb T=U(1)$ and $\rho$ is a character.

\begin{corollary}
\label{dense}
Let the torsion subgroup of an abelian topological group $G$ be
everywhere dense in it. Then every continuous irreducible
representation of $G$ in a Banach space is one-dimensional and indeed a
continuous character.
\end{corollary}

\noindent{\bf Proof.} The restriction of a continuous irreducible
representation $\rho$ of $G$ in a Banach space $E$
to the torsion subgroup $T(G)$ is irreducible because $T(G)$ is
everywhere dense in $G$. Therefore, $E=\C$ and $\rho\vert_{T(G)}$ is
a continuous character by Theorem \ref{t1}. Its extension to a continuous
homomorphism $G\to\Bbb T\subset\mbox{GL}(\C)$ is unique and
must therefore coincide with $\rho$. We conclude that $\rho$ is also a
continuous character, Q.E.D.

\end{section}

\begin{section}{Connected torsion groups}

The following observation was made by A.A. Markov
\cite{M}.

\begin{proposition} 
An abelian torsion group equipped with a connected group
topology is minimally almost periodic.
\label{markov}
\end{proposition}

\noindent{\bf Proof.} The image of $G$ under every continuous
character $\chi\colon G\to\Bbb T$ is connected and consists of
elements of finite order, which means that 
$\chi(G)=\{e_{\Bbb T}\}$ and $\chi$ is trivial.

It follows immediately from Theorem \ref{main} and Proposition \ref{markov}
that all such topological 
groups actually enjoy an (at least, formally) stronger property.

\begin{corollary}
\label{markov2}
An abelian torsion group equipped with a connected group
topology admits no nontrivial
irreducible continuous Banach representations.
\end{corollary}

This result leads to our first class of abelian topological
groups admitting no irreducible representations in Banach spaces.

\begin{examples} {\rm
In \cite{M} Markov had constructed the first ever example of
a connected group topology on an infinite torsion group (whose elements 
have order $2$). He then used the above observation
\ref{markov} to show that $G$ is minimally almost periodic. 
Later Graev \cite{Gr} showed how to produce numerous such examples much more 
easily as varietal free topological 
groups: using the modern terminology \cite{Mo}, every Graev free
topological group, $F_{\frak A_{n}}(X)$ on a nontrivial connected
topological space $X$ formed in the variety $\frak A_{n}$
of all topological groups of finite period $n$ is nontrivial
(moreover, contains $X$ as a topological subspace). 

Still one more way to generate such groups at will
is to embed any torsion group $G$ equipped with the discrete
topology into a path-connected topological group of all $G$-valued
step functions on $[0,1]$ equipped with the
topology of convergence in measure, using the construction of
Hartman and Mycielski \cite{HM}.}
\end{examples}

Our next result follows immediately from Corollary \ref{dense}
and Proposition \ref{markov}.

\begin{corollary}
Let the torsion subgroup of an abelian topological group $G$ be
connected and everywhere dense in $G$. Then $G$ admits no
nontrivial strongly continuous irreducible representations
in Banach spaces.
\label{main}
\end{corollary}
\end{section}

\begin{example} {\rm
Here is an interesting source of rather natural
examples of topological groups possessing the properties stated in
Corollary \ref{main}, which construction
 seems to have never been explored before.
For a pointed topological space $X=(X, \ast)$ denote by $A(X)$ the
{\it Graev free abelian topological group on} $X$,
that is, an abelian topological group algebraically free over 
$X\setminus\{\ast\}$ and
containing $X$ as a closed topological subspace in such a way that
every continuous mapping $f$ from $X$ to an abelian
topological group $G$, sending $\ast$ to the identity,
 gives rise to a unique continuous homomorphism
$\bar f: A(X)\to G$ with $f=\bar f\circ i$. (See e.g.
\cite{Mor}.) Let $L(X)$ be the 
{\it free locally convex space on} $X=(X, \ast)$, that is, a 
locally convex space containing $X\setminus\{\ast\}$
 as a Hamel basis and a closed
topological subspace in such a way that 
every continuous mapping $f$ from $X$ to a
locally convex space $E$, sending $\ast$ to zero,
gives rise to a unique continuous linear operator $\bar f: L(X)\to E$
with $f=\bar f\circ i$. (See \cite{Fl}.)
The identity mapping $id_X:X\to X$ extends to
a canonical continuous homomorphism $i:A(X)\to L(X)$, which is
an embedding of $A(X)$ into the additive 
topological group of 
$L(X)$ as a closed topological subgroup.
(It follows from the corresponding result for the free Markov
abelian topological group and the free Markov locally convex space,
see
\cite{Tk} and \cite{U}.) Denote the factor-group $L(X)/A(X)$
by $\Bbb T(X)$; it is nontrivial whenever $X$ is such. 
(For example, $\Bbb T(\{0,1\})\cong\Bbb T$.)
If $X$ is connected, then the torsion subgroup of
$\Bbb T(X)$ is connected, being algebraically generated by
the image under the quotient homomorphism $L(X)\to\Bbb T(X)$
of the connected set 
$\cup_{n=1}^{\infty}\frac 1 nX$.
}
\end{example}

\begin{section}{Abelian Banach--Lie groups}

\begin{corollary} A minimally almost periodic abelian Banach--Lie group
admits no nontrivial
irreducible continuous Banach representations.
\label{BL}
\end{corollary}

\noindent{\bf Proof.} A minimally almost periodic 
abelian Banach--Lie group is necessarily connected (for otherwise it 
would possess a discrete abelian factor group) and therefore a
quotient group of the additive
group of a Banach space $E$ (the underlying space of the Banach--Lie 
algebra of $G$) modulo a discrete subgroup, $D$. 
Denote by $C$ the union of all
one-parameter subgroups passing through elements of $D\setminus\{0\}$
(that is, by the cone spanned by $D$). 
The closed subgroup of $E$ generated by $C$ is in fact a closed
linear subspace, $F$, and assuming $F\neq E$ leads to a 
contradiction, since the quotient group $E/F$ is a
topological quotient group of $G$ and the additive group of a
non-degenerate Banach space, hence not minimally almost periodic. 
It means that $F=E$ and as a corollary, the image of $C$ 
generates an everywhere dense subgroup of $G$. Since the image of
every one-parameter subgroup of $E$ passing through a point of $D$ is
topologically isomorphic to the group $\Bbb T=U(1)$, the torsion 
subgroup of
$G$ is everywhere dense in $G$. Now Theorem \ref{t1} together
with minimal almost periodicity imply the result.

\begin{examples} {\rm
A large class of examples of abelian
Banach--Lie groups admitting no continuous characters are known,
cf. \cite{Ro}, \cite{DPS}, \cite{Ba}, \cite{Ba2}. Every such group is in fact
{\it monothetic,} that is, contains an everywhere 
dense cyclic subgroup. According to Corollary
\ref{BL}, every such Banach--Lie group admits no irreducible
Banach representations.}
\end{examples}
\end{section}

\begin{section}{Levy groups}

A topological group $G$ is called a {\it Levy group} 
\cite{Gl1, Gl2, Un} if $G$ contains an
increasing chain of compact subgroups $G_i,i\in\N$, having an
everywhere dense union in $G$ and such that whenever 
$\liminf\mu_i(A_i)>0$ for some $A_i\subseteq G_i$, one has
$\lim \mu_i(VA_i\cap G_i)=1$ for every
neighbourhood of identity $V$, 
where $\mu_i$ denotes the normalized Haar measure on $G_i$.

\begin{corollary} An abelian Levy group admits no
irreducible continuous Banach representations.
\label{levy}
\end{corollary}

\noindent{\bf Proof.} Since the torsion subgroup
of every compact abelian group is everywhere dense, 
one concludes that the torsion subgroup of $\cup_{i} G_{i}$ is
everywhere dense in $G$ and
Theorem \ref{t1} applies. The known minimal almost periodicity
of Levy groups accomplishes the proof. 
(Every Levy group $G$ is
{\it extremely amenable,} that is, has a fixed point in every
compact $G$-space, see \cite{Gl2} and also \cite{Un}. This property
obviously implies minimal almost periodicity; whether the converse is
true, is unknown as of August 1997.) 

\begin{example} {\rm
A concrete example of a Levy group $G$
(due to Glasner and, independently, Fursternberg and B. Weiss)
is the abelian group $L_{1}(X,{\Bbb S}^{1})$ of all measurable
$\Bbb S^1$-valued complex functions on a nonatomic Lebesgue measure
space $X$ equipped with the metric
$$d(f,g)=\int_{X}\vert f(x)-g(x)\vert_{\C}~ d\mu(x),$$
where $\Bbb S^{1}$ is identified with the multiplicative subgroup of all
complex numbers of modulus $1$.
As the compact subgroups $G_i$ of $G$, one can choose 
tori of increasing finite dimension formed by
step functions corresponding to a sequence of refining partitions of
$X$. An ingenious argument \cite{Gl2} shows that this group is also
monothetic. }
\end{example}

\end{section}

\begin{section}{Monothetic unitary groups}
There are at least two known
examples of minimally almost periodic monothetic topological groups
admitting a faithful strongly continuous unitary representation in a
Hilbert space. One of them is the Levy group
$L_{1}(X,{\Bbb S}^{1})$, having a faithful continuous unitary representation 
by multiplication operators in the Hilbert space $L_{2}(X)$
\cite{Gl2}. The other known example belongs to Banaszczyk 
(\cite{Ba2}, Th. 5.1), who constructed 
monothetic Banach--Lie groups (modelled over the spaces $l_{p}$) without
continuous characters but admitting a faithful unitary representation.

According to Corollary \ref{levy} and Corollary \ref{BL}, the 
topological groups of the above type admit no nontrivial irreducible Banach
representations. The same property is obviously shared by their
images in the unitary group $U(\cal H)$ 
equipped with the strong operator topology. We have established the
following.

\begin{corollary}
\label{comb}
There exist monothetic topological subgroups $G$ of the 
unitary group $U(\cal H)$ of the separable Hilbert space
equipped with the strong operator topology that admit no
nontrivial irreducible strongly continuous Banach representations.
\end{corollary}

\end{section}

\begin{section}{Acknowledgements}
I am grateful to the Editors of
the book `Abelian Groups, Module Theory and Toplogical 
Algebra' Professors Dikran Dikranjan and Luigi Salce for inviting 
me to submit a paper, and to Professor Wojciech Banaszczyk 
for a useful remark.
\end{section}

\noindent
School of Mathematical and Computing Sciences,
Victoria University of Wellington, P.O. Box 600,
Wellington, New Zealand\\
e-mail address: vladimir.pestov@vuw.ac.nz\\

\end{document}